# The Location-Allocation Model for Multi-Classification-Yard Location Problem in a Railway Network

Boliang Lin[*]

* School of Traffic and Transportation, Beijing Jiaotong University, Beijing 100044, bllin@bjtu.edu.cn;

**Abstract:** Classification yards are crucial nodes of railway freight transportation network, which plays a vital role in car flow reclassification and new train formation. Generally, a modern yard covers an expanse of several square kilometers and costs billions of yuan, i.e., hundreds of millions of dollars. The determination of location and size of classification yards, which is a location-allocation problem with railway characteristics, is not only related to building or improving cost, but also involved with train connecting service (TCS) plan. This paper proposed a bi-level programming model for this problem. The upper-level is intended to find an optimal building or improving strategy for potential nodes, and the lower-level aims to obtain a least costly TCS plan considering reclassification cost and accumulation delay, when the building or improvement plan is given by the upper-level. The model is constrained by capital budget, classification capacity, the number of available tracks, etc.

## 1. Introduction

Classification yards are usually referred as the hubs for freight train formation, handling a great many of freight trains. The spatial configuration of yards and classification workload distribution among them play a vital role in increasing transportation efficiency and benefits, improving the utilization of rolling stocks (locomotives and railcars), and promoting the development of railway industry. Thus, the determination of location and size of yards is a crucial problem in network design. Generally, a modern classification yard costs hundreds of millions of dollars and occupies a land of several square kilometers. For instance, the length of Maschen marshalling yard in Germany reaches 7000 m, the width of which is up to 700 m, covering 2.8 km$^2$. The world's largest classification yard, Bailey Yard in the United

States, has a length of about 13000 m, and a maximum width of 3200 m, taking up an area of 11.5 km$^2$. In China, Wuhan North Railway marshalling station has a length of over 5000 m and a width of nearly 1000 m, with a land occupation of some 4.5 km$^2$. Zhengzhou North Classification Yard, once the largest one in Asia, covers a total expanse of 5.3 km$^2$, and is over 6000 m in length and 800 m in width. To avoid huge waste of resources due to unreasonable building and improving, scientific location strategy has already been the important basis in decision making for departments concerned. Hence, the research on classification yard location problem is of great value.

The cost of building or improving a marshalling station constitutes a significant portion of capital investment. And the location, automaticity and classification ability directly affects the performance and operation cost of the yard. Furthermore, the spatial configuration of yards significantly influences the routing of traffic flows over the whole network. Therefore, a reasonable location and appropriate size is a guarantee for considerable profits. Given the highly-nonlinear interrelation among yards, the freight train formation plan needs to be taken into account in investment analysis, which should be carried out from the perspective of railway network rather than focusing on a certain yard. As the budget is limited, it has to be ensured that the new-building or improvement plan to be carried out can best meet the goals of investors. From the micro-economic perspective, yards should handle their workloads in minimum unit cost. Note that, in comparison with the unit cost of yard having optimal workload-capacity ratio, the classification cost per car is larger in the yard whose available capacity is much larger than workloads, which is a huge waste of classification resources. From the macro-economic perspective, capital investment is expected to improve the routing of railcars over the network, i.e., the vast majority of railcars can be reclassified in optimal yards respectively, and the capacity of them is rationally utilized.

Theoretically, the number of combination strategies for building or improving classification yards grows exponentially with the number of candidate locations. Let us assume that there are ten potential locations, each having three new-building or improvement strategies. For this case, the total number of combination schemes reaches 59049. Furthermore, each scheme corresponds to a freight train connecting service problem, resulting in huge computational cost. Therefore, the yard location problem is of high complexity, featuring profound theoretical depth and great application value.

To summarize, quantitatively analyzing the multi-classification-yard location problem, from the perspective of capital investment and train formation cost, has already become a problem to be solved urgently.

## 2. Literature Review

Many researchers have studied in-depth the LAP. Work on LAP was initiated by Cooper [1] who presented exact extremal equations and a heuristic method for solving certain classes of LAP. Bongartz et al. [2] proposed a solution method which relaxed the 0-1 constraints on the allocations and solved for both the locations and allocations simultaneously. Eben-Chaime et al. [3] studied the capacitated LAP on a line, established mathematical optimization models and discussed their properties and complexity.

Research on location problem of classification yard started pretty early which could be traced back to the 1950s. Mansfield and Wein [4] put forward the first location model of classification yard in 1958. The model is constructed to aid a railroad management in choosing among alternative locations when newer facilities were going to be installed. However, the background that a large amount of classification yards operated and the freight trains were reclassified almost at every yard they passed through, was completely different from current situation. Assad [5] proposed the general principles of yard location, but did not give the method of determining the yard quantity and scale. Maji and Jha [6] constructed a location model of classification yard aiming at minimizing the sum of fixed cost and variable cost. The reliability and effectiveness of the model were tested on a railway network composed of 16 stations. Lee et al. [7] developed an optimal marshalling yard location model considering economies of scale due to the consolidation of flows. The model was tested on a regional railway network consisting of 25 nodes in South Korea, and the result was compared with that of O'Kelly, Racunica and Wynter, respectively. Lin et al. [8] proposed a multi-project decision model to determine the location and size of new classification yards as well as the technical upgrading plan of existing yards. The model took into account the capital investment, reclassification cost and accumulation cost, which differentiated it from other models for classification yard location problem. Yan et al. [9] established a typical station type decision model considering wagon flow structure. Li et al. [10] constructed a 0-1 programming model for train formation plan and classification yard location. Likewise, Geng [11] developed an integrated optimization model for freight train connecting service plan and reclassification load distribution. In addition, the bi-level programming model has also been extensively studied by many experts, and has achieved good results [12-14].

The paper is organized as follows: Section 1 describes the importance and complexity of the multi-classification-yard location problem in detail. Section 2 provides a brief survey of the literature devoted to the location-allocation problem and classification yard location problem. In Section 3 we establish a bi-level programming model in order to obtain a good investment strategy. And conclusions are presented in Section 4.

## 3. Mathematical Model

In this section, we first present detailed descriptions of the notations used in this paper. After that, a bi-level programming model is proposed for the multi-classification-yard location problem. Some discussions with respect to the model are also presented in this section.

### 3.1 Notations

The notations used in this paper are listed in Table 1.

**Table 1** Notations used in this paper.

| Sets | Descriptions |
|---|---|
| $V$ | The set of all nodes in a rail network; |
| $V^{Potential}$ | The set of all potential nodes where yards are built or improved, including all candidate new yard locations and existing yards that may be improved. Obviously, $V^{Potential} \subseteq V$; |
| $V^{Original}$ | The set of original yards in a rail network; |
| $P(k)$ | The set of investment plans for yards building or improvement at node $k$, excluding the plan of no investment; |
| $\rho_{ij}$ | The set of yards through which a flow from $i$ to $j$ pass on its itinerary, excluding yard $i$ and yard $j$. |
| **Parameters** | **Definitions** |
| $I_k^p$ | The cost of investment plan $p$ at node $k$; |
| $B$ | Budget of investment; |
| $T_k^p$ | The lifetime of node $k$ when choosing plan $p$; |
| $\gamma$ | The discount rate of capital investment; |
| $\alpha$ | The coefficient of converting car-hour cost into economic cost; |
| $c_i$ | The accumulation parameter of yard $i$; |
| $m_{ij}$ | The size of train dispatched from $i$ to $j$; |
| $F_k$ | The workload of node $k$, i.e., the number of cars reclassified at node $k$; |
| $\tau_k$ | The original classification cost per railcar at node $k$ before improvement; |
| $\tau_k^p$ | The classification cost per railcar at node $k$ after building or improvement based on plan $p$; |
| $C_k^{\text{Total}}$ | The classification capacity of node $k$ before building or improvement; |
| $C_k^{\text{Local}}$ | The classification capacity reserved for local trains at node $k$ before building or improvement; |
| $\Delta C_k^p$ | The increase of classification capacity at node $k$ after building or improvement based on plan $p$; |
| $L_i^{\text{Total}}$ | The number of classification tracks at node $i$ before building or improvement; |
| $L_i^{\text{Local}}$ | The classification tracks reserved for local trains at node $i$ before building or improvement; |
| $\Delta L_i^p$ | The increase of classification tracks at node $i$ after building or improvement based on plan $p$; |
| $N_{ij}$ | The traffic demand which origins at node $i$ and is destined to node $j$; |
| **Decision variables** | **Definitions** |
| $y_k^p$ | Investment variable; it takes value one if plan $p$ is selected for node $k$, and zero otherwise. |
| $x_{ij}$ | Train variable; its value is one if the train service $i \to j$ is provided, and zero otherwise. |
| $x_{ij}^k$ | Car flow variable; it takes value one if cars whose destination is $j$ are consolidated into train service $i \to k$ at node $i$. Otherwise, it is zero. |

## 3.2 Model Descriptions

For a candidate yard location or an existing yard that might be improved, $k \in V^{\text{Potential}}$, the set of its investment plans is denoted as $P(k)$. For instance, if node $k$ is not a classification yard at present, it may have three new-building plans for selection, including single directional longitudinal-type marshalling station with three yards in three stages, single directional combination-type marshalling station with four yards in two stages, and double directional longitudinal-type marshalling station with six yards in three stages. Conversely, if node $k$ is already a single directional longitudinal-type marshalling station with three yards in three stages, its investment plan might be improved into a double directional longitudinal-type marshalling station with six yards in three stages. As the investment plans of nodes are not necessarily selected, a constraint should be considered:

$$\sum_{p \in P(k)} y_k^p \leq 1 \quad \forall k \in V^{\text{Potential}} \quad (1)$$

To convert the inequality constraint mentioned above into a equality constraint, the plan of no investment can be viewed as a special plan "0" whose corresponding cost is $I_k^0 = 0$, and the decision variable is denoted as $y_k^0$. Thus, the constraint (1) can be modified into:

$$y_k^0 + \sum_{p \in P(k)} y_k^p = 1 \quad \forall k \in V^{\text{Potential}} \quad (2)$$

Note that the set of decision variables of a certain node can be denoted by $Y_k = (y_k^0, y_k^1, \cdots, y_k^{|P(k)|})$. For this case, $Y = (Y_1, Y_2, \cdots, Y_k \ldots, Y_{|V^{\text{Potential}}|})$ constitutes a decision space for a group of candidate nodes.

If we choose plan $p$ for a certain node $k$, i.e., $y_k^p = 1$, then its capital investment is $I_k^p$, and the lifetime is set to $T_k^p$. Without loss of generality, it is assumed that, at the end of its lifetime, the salvage value of a yard will decrease to zero (in fact, once the yard is built, it may exist forever unless been closed, as new equipments are acquired when needed). Capital recovery factor is introduced to annualize the project investment $I_k^p$, which can be described by $\dfrac{\gamma(1+\gamma)^{T_k^p}}{(1+\gamma)^{T_k^p}-1} I_k^p$.

Given a set of investment plans, the daily operation cost for all shipments, including accumulation delay and reclassification delay, is denoted as $Z(Y)$. Therefore, an upper-level formulation can be constructed to describe the multi-classification-yard location problem:

**Upper-level program:**

$$\min \sum_{k \in V^{\text{Potential}}} \sum_{p \in P(k)} \frac{\gamma(1+\gamma)^{T_k^p}}{(1+\gamma)^{T_k^p}-1} I_k^p y_k^p + 365 \times \alpha Z(Y) \tag{3}$$

s.t.

$$y_k^0 + \sum_{p \in P(k)} y_k^p = 1 \quad \forall k \in V^{\text{Potential}} \tag{4}$$

$$\sum_{k \in V^{\text{Potential}}} \sum_{p \in P(k)} I_k^p y_k^p \leq B \tag{5}$$

$$y_k^0, y_k^p \in \{1,0\} \quad \forall k \in V^{\text{Potential}}, p \in P(k) \tag{6}$$

The first term of the objective function is the annualized project investment for all potential nodes. The second term is the annual operation cost. The constraint (5) ensures that the capital investment of all potential nodes would not exceed the budget. Apparently, without consideration of the operation cost reduction due to yards building or improvement, the upper-level reaches its optimality when the investment for all potential nodes is zero, i.e., $y_k^0 = 1$, $I_k^0 = 0$ $(\forall k \in V^{\text{Potential}})$. For example, some shipments that should be reclassified at yard $k$ (according to the optimal train connecting service plan), have to be reclassified at other yards due to the capacity constraint of yard $k$. In this case, the operation cost will definitely increase. In other words, no project investment is the best plan when not considering the cost savings associated with operation efficiency improvement. Indeed, there was no hump yard at the initial stage of railway industry. With the increase of traffic volume and classification workload, people began to invest huge amount of money to build hump yards.

Typically, establishing or improving a yard will not only increase classification capacity $\Delta C_k^P$ and the number of tracks $\Delta L_i^p$ (more blocks can be built simultaneously), but also raise operation efficiency (due to the improvement of automaticity) and reduce classification cost $\tau_k^p$. Therefore, the determination of train connecting service and distribution of classification workload (the minimization of $Z(Y)$), in the case of given location strategy, can be referred as an allocation problem. If a train connecting service is provided between two yards, $x_{ij} = 1$, then an accumulation delay $c_i m_{ij}$ will be incurred at the origin yard. For instance, if

$c_i = 12$, $m_{ij} = 50$, the accumulation delay will be 600 car-hours. In fact, the practical accumulation parameter $c_i$ is generally less than 12 due to accumulation disruption. It is clear that the fixed charge of providing a train service is associated with the size of train rather than traffic volume (details can be referred to Lin et al., 2012). When a shipment is reclassified at yard $k$ (carried by a inbound train, then reclassified and assembled into an outbound train), each car of the shipment will incur a reclassification cost $\tau_k$ (including time cost, labor expense and fuel consumption). Likewise, if plan $p$ is selected for yard $k$ in the upper-level program ($y_k^p = 1$), the new reclassification cost will be $\tau_k^p$. For the sake of simplicity and clarity, we measure the cost in car-hour consumption.

**Lower-level program:**

$$\min \sum_{i \neq j \in V^{\text{Original}} \bigcup V^{\text{Potential}}} c_i m_{ij} x_{ij} + \sum_{k \in V^{\text{Original}} - V^{\text{Potential}} \bigcap V^{\text{Original}}} F_k \tau_k + \sum_{k \in V^{\text{Potential}}} F_k \left( \sum_{p \in P(k)} \tau_k^p y_k^p + \tau_k y_k^0 \right) \quad (7)$$

s.t.

$$x_{ij} + \sum_{k \in \rho_{ij}} x_{ij}^k = 1 \qquad \forall\, i \neq j \in V^{\text{Original}} \bigcup V^{\text{Potential}} \tag{8}$$

$$F_k \leq C_k^{\text{Total}} - C_k^{\text{Local}} \qquad \forall\, k \in V^{\text{Original}} - V^{\text{Original}} \bigcap V^{\text{Potential}} \tag{9}$$

$$F_k \leq C_k^{\text{Total}} - C_k^{\text{Local}} + \sum_{p \in P(k)} y_k^p \Delta C_k^p \qquad \forall k \in V^{\text{Potential}} \tag{10}$$

$$\sum_{j \in V^{\text{Original}} \bigcup V^{\text{Potential}}} \varphi(D_{ij}) \leq L_i^{\text{Total}} - L_i^{\text{Local}} \qquad \forall\, i \in V^{\text{Original}} - V^{\text{Original}} \bigcap V^{\text{Potential}} \tag{11}$$

$$\sum_{j \in V^{\text{Original}} \bigcup V^{\text{Potential}}} \varphi(D_{ij}) \leq L_i^{\text{Total}} - L_i^{\text{Local}} + \sum_{p \in P(i)} y_i^p \Delta L_i^p \qquad \forall i \in V^{\text{Potential}} \tag{12}$$

$$x_{ij}, x_{ij}^k \in \{1, 0\} \qquad \forall i \neq j, k \in V^{\text{Original}} \bigcup V^{\text{Potential}} \tag{13}$$

The objective function of the lower-level consists of three terms. The first term is the total accumulation delay cost of all train services. The second term is the classification cost of yards not included in $V^{\text{Potential}}$. While the third term is the classification cost of nodes included in $V^{\text{Potential}}$. $F_k$ is the classification workload at node $k$ which can be described by:

$$F_k = \sum_i \sum_j f_{ij} x_{ij}^k \qquad \forall i \neq j, k \in V^{\text{Original}} \bigcup V^{\text{Potential}} \tag{14}$$

where $f_{ij}$ is the car flow from yard $i$ to yard $j$, consisting of the original demand $N_{ij}$ and the reclassified cars from other yards. It can be expressed by:

$$f_{ij} = N_{ij} + \sum_{s \in V^{\text{Original}} \bigcup V^{\text{Potential}}} f_{sj} x_{sj}^{i} \quad \forall i \neq j \in V^{\text{Original}} \bigcup V^{\text{Potential}} \tag{15}$$

The constraint (8) guarantees that a car flow can either be directly shipped to the destination or classified at more than one intermediate yard on its itinerary. For the yards not included in $V^{\text{Potential}}$, the constraint (9) ensures that the workload would not exceed their capacity, and the constraint (11) guarantees that the occupied tracks is less than the number of available classification tracks. Similarly, for the nodes included in $V^{\text{Potential}}$, the constraint (10) ensures that the capacity of these nodes is enough to handle the workload, and constraint (12) guarantees that there are enough available classification tracks for storing outbound trains.

Note that, $D_{ij}$ is defined as the service flow from yard $i$ to $j$, which means the number of cars shipped by train service $i \rightarrow j$, and it can be described as follows:

$$D_{ij} = f_{ij} x_{ij} + \sum_{t \in V^{\text{Original}} \bigcup V^{\text{Potential}}} f_{it} x_{it}^{j} \quad \forall i \neq j \in V^{\text{Original}} \bigcup V^{\text{Potential}} \tag{16}$$

The relation between service flow $D_{ij}$ and the number of tracks needed, namely the track demand function $\varphi(D_{ij})$ can be expressed by:

$$\varphi(D_{ij}) = D_{ij} / 200 \tag{17}$$

or described by:

$$\varphi(D_{ij}) = \begin{cases} 1 & 0 < D_{ij} \leq a_1 \\ 2 & a_1 < D_{ij} \leq a_2 \\ \vdots & \vdots \\ n & a_{n-1} < D_{ij} \leq a_n \end{cases} \tag{18}$$

where the values of $a_1, a_2, \cdots, a_n$ should logically depend on the yard equipment and labor resources. In China railroad system, we usually set $a_1 = 200, a_2 = 400, \cdots, a_n = 200n$. The equation (17) is adapted in Lin et al. (2012). It should be noted that the calculation result based on equation (17) might have decimal part, which will not occur in formulation (18). In fact, a classification track, in practice, cannot simultaneously store railcars that are assigned to different outbound trains, even though the capacity of the track is enough to handle these railcars. Thus, it seems that formulation (18) is more practical and reasonable.

## 4. Conclusions

In this paper, we model the multi-classification-yard location problem as a location-allocation problem with railway characteristics, and proposed a bi-level programming model constrained by budget, classification capacity and number of available tracks. The upper-level is intended to find an optimal set of investment plans for all potential nodes, and the lower-level aims to obtain a least costly train connecting service plan considering reclassification cost and accumulation delay, in the context of given investment plans. Considering the highly-nonlinear interrelation among yards, the investment plan should be analyzed from the perspective of railway network as a whole rather than focusing on a certain yard. For simplicity, the plan of no investment is incorporated in the model as a special plan. And the capital investment for all potential nodes is annualized by capital recovery factor.

In the long term, researcher can focus on the multistage investment problem for classification yards. We identify this a promising area for future research.

## 5. References

[1] Cooper L. Location-Allocation Problems [J]. Operations Research, 1963, 11(3):331-343.

[2] Bongartz I, Calamai P H, Conn A R. A projection method for, 1 p, norm location-allocation problems [J]. Mathematical Programming, 1994, 66(1):283-312.

[3] Eben-Chaime M, Mehrez A, Markovich G. Capacitated location-allocation problems on a line [J]. Computers & Operations Research, 2002, 29(5):459-470.

[4] Mansfield E, Wein H H. A Model for the Location of a Railroad Classification Yard[J]. Management Science, 1958, 4(3):292-313.

[5] Assad A. Models for Rail Transportation [J]. Trans. ReS., 1980, 14A: 205-220.

[6] Avijit Maji, Manoj K. Jha. Railroad Yard Location Optimization Using A Genetic Algorithm [J]. Environmental Science and Sustainability.

[7] Lee J S, Kim D K, Chon K S. Design of Optimal Marshalling Yard Location Model Considering Rail Freight Hun Network Properties [J]. Journal of the

Eastern Asia Society for Transportation Studies, 2008, 7:1031-1045.

[8] Lin B.L., Xu Z.Y., Huang M., Ji J.L., Guo P.W. An optimization approach to railroad classification yard location and installment size decision with budget constraint [J]. Journal of the China Railway Society, 2002,24(3),5-8.

[9] Yan H.X., Lin B.L., Zheng J.R., Zhang Y.H. Comparative study on optimal type selection of marshalling station between unidirectional and bidirectional [J]. Journal of transportation systems engineering and information technology, 2007, 7(1), 124-131.

[10] Li H.D., He S.W., Song R., Ji L.J., Shen Y.S. Synthetic optimization of train formation plan and layout of technical service station [J]. Journal of Beijing Jiaotong University, 2010, 34(6), 30-34+39.

[11] Geng L.Q. Optimization model of freight car marshalling scheme and load dividing and combining in marshalling station [J]. Railway transport and economy, 2011, 33(6), 59-63+85.

[12] Shi F., Fang Q.G., Li X.H., Mo H.H., Huang Y.L. A bi-level programming optimization method for the layout of technical service station [J]. Journal of the China Railway Society, 2003, 25(2), 1-4.

[13] Yin Y., Peng Q.Y. Bi-level programming method of railway marshalling station layout [J]. Railway transport and economy, 33(8), 25-29.

[14] Li S. Combinatorial research of rail marshalling yard layout and railcar blocking problems [D]. Beijing: Beijing Jiaotong University, 2014.

[15] O'Kelly M E．The location of interacting hub facilities[J]．Transportation Science，1 986，20(2)：92-1 06．

[16] O'Kelly M E．A quadratic integer program for the location of interacting hub facilities[J]．European Journal of Operational Research，1987，32(3)：393-404．

[17] Campbell J F, Ernst A T, Krishnamoorthy M. Hub Arc Location Problems: Part I-Introduction and Results [J]. Management Science, 2005, 51(10):1540-1555.

[18] Campbell J F, Ernst A T, Krishnamoorthy M. Hub Arc Location Problems: Part II: Formulations and Optimal Algorithms [J]. Management Science, 2005, 51(10):1556-1571.

[19] Racunica I, Wynter L. Optimal location of intermodal freight hubs [J]. Transportation Research Part B Methodological, 2005, 39(5):453-477.

[20] Arnold P, Peeters D, Thomas I. Modelling a rail/road intermodal transportation

system [J]. Transportation Research Part E Logistics & Transportation Review, 2004, 40(3):255-270.

[21] Ebery J, Krishnamoorthy M, Ernst A, et al. The capacitated multiple allocation hub location problem: Formulations and algorithms [J]. European Journal of Operational Research, 2000, 120(3):614–631.

[22] Khalafallah A, El-Rayes K. Automated multi-objective optimization system for airport site layouts[J]. Automation in Construction, 2011, 20(4): 313-320.

[23] Oktal H, Ozger A. Hub location in air cargo transportation: A case study [J]. Journal of Air Transport Management, 2013, 27: 1-4.

[24] Figueiredo R, O'Kelly M E, Pizzolato N D. A two-stage hub location method for air transportation in Brazil [J]. International Transactions in Operational Research, 2014, 21(2): 275-289.

[25] Yang T H. Airline network design problem with different airport capacity constraints [J]. Transportmetrica, 2008, 4(1): 33-49.

[26] Yang Z, Yu S, Notteboom T. Airport location in multiple airport regions (MARs): The role of land and airside accessibility [J]. Journal of Transport Geography, 2016, 52: 98-110.